\documentclass[reqno,11pt]{amsart}

\usepackage{a4,amsmath,amssymb,amsthm,eucal}
\usepackage{graphicx}
\usepackage{pdfsync}
\usepackage{epstopdf}
\def\R{\mathbb{R}}

\def\A{{\mathcal A}}
\def\B{{\mathcal B}}
\def\C{{\mathcal C}}

\def\W{{\mathcal W}}

\def\argmin{\mathrm{argmin}\,}
\def\argmax{\mathrm{argmax}\,}

\def\lip{\mathrm{Lip}}
\def\eps{\varepsilon}

\newtheorem{theorem}{Theorem}[section]
\newtheorem{definition}[theorem]{Definition}

\newtheorem{lemma}[theorem]{Lemma}
\newtheorem{proposition}[theorem]{Proposition}

\theoremstyle{remark}
\newtheorem{remark}[theorem]{Remark}
\newtheorem{example}[theorem]{Example}

\newcommand\pref[1]{(\ref{#1})}
\newcommand\tilp{\widetilde{p}}
\newcommand\tilu{\widetilde{u}}

\newcommand\Omb{\overline{\Omega}}

\numberwithin{equation}{section}

\author{G. Buttazzo\and G. Carlier}
\thanks{\it This research has been conceived during a visit of the first author to CEREMADE of Universit\'e de Paris Dauphine; he wishes to thank this institution for the warm and friendly atmosphere provided during all the visit.}
\title[Optimal spatial pricing strategies with transportation costs]{Optimal spatial pricing strategies with transportation costs}

\begin{document}

\begin{abstract}
We consider an optimization problem in a given region $Q$ where an agent has to decide the price $p(x)$ of a product for every $x\in Q$. The customers know the pricing pattern $p$ and may shop at any place $y$, paying the cost $p(y)$ and additionally a transportation cost $c(x,y)$ for a given transportation cost function $c$. We will study two models: the first one where the agent operates everywhere on $Q$ and a second one where the agent operates only in a subregion. For both models we discuss the mathematical framework and we obtain an existence result for a pricing strategy which maximizes the total profit of the agent. We also present some particular cases where more detailed computations can be made, as the case of concave costs, the case of quadratic cost, and the onedimensional case. Finally we discuss possible extensions and developments, as for instance the case of Nash equilibria when more agents operate on the same market.
\end{abstract}

\maketitle

\section{Introduction}\label{sec1}

In the present paper we consider a model where in a prescribed region $Q$ of the Euclidean space $\R^d$ an agent (a central government or a commercial company) has the possibility to decide the price of a certain product; this price $p(x)$ may vary at each point $x\in Q$ and the customers density $f(x)$ is assumed to be completely known.

We assume that all the customers buy the same quantity of the product; on the counterpart, a customer living at the point $x\in Q$ knows the pricing function $p$ everywhere and may decide to buy the product where he lives, then paying a cost $p(x)$, or in another place $y$, then paying the cost $p(y)$ and additionally a transportation cost $c(x,y)$ for a given transportation cost function $c$.

The individual strategy of each customer is then to solve the minimization problem
\begin{equation}\label{pbT}
\min_{y\in Q}\Big\{c(x,y)+p(y)\Big\}.
\end{equation}

Of particular importance to our problem is the  (set-valued) map  $T_p: Q \to Q$ which associates to every customer living at the point $x$ all the locations where it is optimal to purchase the good. Given the price pattern $p$, $T_p$ is then defined by
\begin{equation}\label{tp}
T_p(x):=\argmin_{y\in Q}\Big\{c(x,y)+p(y)\Big\},\; \forall x\in Q.
\end{equation}

Without any other constraint, due to the fact that the customers have to buy the product (for instance  gasoline,  food,  a medical product or  cigarettes), the pricing strategy for the agent in order to maximize the total income would simply be increasing everywhere the function $p$ more and more. To avoid this trivial strategy we assume that on the region $Q$ some kind of regulations are present, and we study the optimization problems the agent has to solve in order to maximize its total profit. We will study two models according to two different price constraints. We also assume that the supply is unconstrained at any location of the region $Q$, which means that whatever the total demand for the product is at a given location, it can be supplied by the agent to the customers. 

\subsection{The agent operates everywhere}\label{modone}

The simplest situation we consider is when the price $p(x)$ is constrained to remain below a fixed bound $p_0(x)$ everywhere on $Q$, due for instance to some regulatory policy. The only assumption we make is that $p_0$ is a {\it proper} nonnegative function, intending that in the region where $p_0=+\infty$ no restrictions on $p(x)$ are imposed. The goal of the agent is to maximize its total income that, with the notation introduced in \eqref{pbT} and \eqref{tp}, can be written as
\begin{equation}\label{incomeone}
F(p,T):=\int_Q p(Tx)\,df(x)
\end{equation}
under the constraint (state equation) that $Tx\in T_p(x)$ i.e. that $T$ is compatible with the customer's individual minimization problem \eqref{pbT}. One may therefore see the previous program as a nonstandard optimal control problem where $p$ is the control and $T$ the state variable. Let us mention that problems with a similar structure naturally arise in the so-called {\it principal-agent} problem in Economics (see for instance Rochet and Chon\'e \cite{rochetchone} and the references therein).

\subsection{The agent operates in a subregion}\label{modtwo}

We consider a second model of pricing strategy: we suppose that in $Q$ there is a given subregion $Q_0$ where the price $p(x)$ is fixed as a function $p_0(x)$ that the agent cannot control. This is for instance the case of another country if the agent represents a central government, or of a region where for some social reasons that the agent cannot modify, the prices of the product are fixed.

Whenever $T_p(x)\subset Q_0$, then the agent makes no benefit from customers living at $x$. In fact the total profit of the agent is given by
\begin{equation}\label{income}
\Pi(p,T):=\int_{T^{-1}(Q\setminus Q_0)} p(Tx)\,df(x).
\end{equation}
under the constraint (state equation) that $Tx\in T_p(x)$. Note that in  formula \eqref{income} giving the total profit, the integration is now performed only on the set of customers that do shop in the region controlled by the agent and not in the fixed-price region $Q_0$. The problem we are interested in reads again as the maximization of the functional $\Pi(p,T)$ among the admissible choices of state and control variables.

\bigskip

For both  models above we discuss the mathematical framework which enables us to obtain an existence result for an optimal pricing strategy and we present some particular cases where more detailed computations can be made, as the case of concave costs, the case of a quadratic cost, and the onedimensional case.

The last section contains some discussions about possible extensions and developments, as for instance the case of Nash equilibria when more agents operate on the same market.

\section{Problem formulation in the first case}\label{firstpb}

In what follows, $Q$ will be some compact metric space (the economic region), and $p_0:Q\to[0,+\infty]$ a nonnegative proper function, i.e. we assume that $p_0$ is not $+\infty$ everywhere on $Q$. We are also given a transportation cost function $c$ assumed to be continuous and nonnegative on $Q\times Q$ and such that $c(x,x)=0$ for all $x\in Q$. Finally, $f$ is a nonnegative Radon measure on $Q$ that models the distribution of customers in $Q$.

The unknown of the problem is the pricing pattern $p$ that varies in the class
$$\A:=\big\{p:Q\to\overline\R\ :\ p\le p_0\mbox{ on }Q,\ p\mbox{ l.s.c. on }Q\big\}.$$
Once a price $p$ has been chosen by the agent, consumers living at any point $x\in Q$ purchase the good so as to minimize their total expenditure which is given by price plus commuting cost. This leads to the following definitions:
\begin{equation}\label{defvptp}
\left\{\hskip-0.4truecm\begin{array}{ll}
&v_p(x):=\min_{y\in Q}\big\{c(x,y)+p(y)\big\},\\
&T_p(x):=\big\{y\in Q\ :\ c(x,y)+p(y)=v_p(x)\big\}.
\end{array}\right.
\end{equation}
By our l.s.c. and compactness assumptions, $T_p(x)$ is a nonempty compact subset of $Q$ (but $T_p$ is not single-valued in general); moreover the graph of $T_p$ is compact as the $\argmin$ of some l.s.c. function on $Q\times Q$. Note that $T_p(x)$ is the set of locations where consumers living at $x$ rationally choose to purchase the good. It is possible however that, for some customers, the optimal total cost $v_p(x)$ is reached at more than one point $y\in Q$; in this case we assume the {\it tie-breaking rule} that the consumers living at $x$ choose to go to a transportation-minimizing (or equivalently price-maximizing) location $y\in T_p(x)$:
$$y\in\argmin_{T_p(x)}c(x,\cdot)=\argmax_{T_p(x)}p(\cdot)$$ 
(notice that every $y\in T_p(x)$ yields the same minimal total expenditure to the customer living at $x$). 

With the previous notation, the optimal pricing problem amounts to
\begin{equation}\label{maxp}
\max\Big\{F(p)\ :\ p\in\A\Big\}
\end{equation}
where $F$ is the functional
\begin{equation}\label{functf}
F(p)=\int_Q\Big(\max_{y\in T_p(x)} p(y)\Big)\,df(x).
\end{equation}
By the definition of $T_p$ and $v_p$, one has $v_p(x)=c(x,y)+p(y)$ for all $y\in T_p(x)$, hence the profit functional can be rewritten as: 
$$F(p)=\int_Q\Big(v_p(x)-\min_{y\in T_p(x)} c(x,y)\Big)\,df(x).$$

In order to obtain the existence of a solution to the optimization problem \eqref{maxp} we reformulate the problem by using the variable $v=v_p$ instead of $p$; the advantage is that $v$ is searched among $c$-concave functions, while $p$ does not have special properties, and this will enable us to obtain the extra compactness necessary to prove the existence result.

\begin{definition}\label{cconc}
A function $v:Q\to\R$ is called $c$-concave if there exists a function $u:Q\to\R$ (that without loss of generality can be assumed upper semicontinuous) such that 
\begin{equation}\label{cconcave}
v(x)=\inf\big\{c(x,y)-u(y)\ :\ y\in Q\big\}.
\end{equation}
For every $c$-concave function $v$ the $c$-transform $v^c$ is defined by
$$v^c(y)=\inf\big\{c(x,y)-v(x)\ :\ x\in Q\big\}$$
and the $c$-superdifferential $\partial^c v(x)$ is given by
$$\partial^c v(x)=\big\{y\in Q\ :\ v(x)+v^c(y)=c(x,y)\big\}.$$
\end{definition}

The previous definition expresses that $c$-concave functions are functions that can be written as pointwise infima of functions of the form $x\mapsto c(x,y)-u(y)$ for some $u:Q\to\R$. The analogy with concave functions (as infima of affine functions) and the parallel between the $c$-transform and the more familiar Legendre-Fenchel transform  then should be clear to the reader. In a similar way, the notion of $c$-superdifferential generalizes the notion of superdifferential for a concave function, and one can actually characterize  $c$-superdifferentials in terms of  the so-called $c$-cyclical monotonicity property that is analogous to the usual cyclical monotonicity. Let us remark as a first example that if $c$ is a distance, then $c$-concave functions are exactly $1$-Lipschitz functions, and in this case one can take $u=-v$ in \pref{cconcave}.  The case of strictly convex costs, and in particular the quadratic cost, will be treated in subsection \ref{strconvex}.

\begin{lemma}\label{unifcont}
Every $c$-concave function is uniformly continuous and its continuity modulus is bounded by the continuity modulus of the cost function $c$ on $Q\times Q$.
\end{lemma}

\begin{proof}
Take a $c$-concave function $v$ and two points $x_1,x_2\in Q$. By the definition of $c$-concavity, for a suitable upper semicontinuous function $u$ we have
$$v(x_2)=\min\big\{c(x_2,y)-u(y)\ :\ y\in Q\big\}=c(x_2,y_2)-u(y_2)$$
where $y_2$ is a suitable point in $Q$. Then we have
$$v(x_1)\le c(x_1,y_2)-u(y_2)=v(x_2)+c(x_1,y_2)-c(x_2,y_2).$$
Interchanging the role of $x_1$ and $x_2$ we deduce
$$|v(x_1)-v(x_2)|\le|c(x_1,y_2)-c(x_2,y_2)|$$
which concludes the proof.
\end{proof}

\begin{lemma}\label{vc}
If $(v_n)$ is a sequence of $c$-concave functions converging uniformly to $v$, then $v_n^c$ converge uniformly to $v^c$ and $v$ is $c$-concave. As a consequence, for every $x\in Q$ we have
$$\min_{y\in\partial^cv(x)}c(x,y)\le\liminf_n\Big(\min_{y\in\partial^cv_n(x)}c(x,y)\Big).$$
\end{lemma}

\begin{proof}
Since $v_n^c$ are $c$-concave, by Lemma \ref{unifcont} it is enough to show that $v_n^c$ converge to $v^c$ pointwise on $Q$. Fix $y\in Q$ and let $x_n\in Q$ be such that
$$v_n^c(y)=c(x_n,y)-v_n(x_n).$$
Since $Q$ is compact, a subsequence of $(x_n)$ converges to some $x\in Q$, so that
$$\liminf_n v_n^c(y)=c(x,y)-v(x)\ge v^c(y).$$
Vice versa, if $x\in Q$ is such that $v^c(y)=c(x,y)-v(x)$, we have $v_n^c(y)\le c(x,y)-v_n(x)$, so that
$$\limsup_n v_n^c(y)\le c(x,y)-v(x)=v^c(y).$$
The fact that $v$ is $c$-concave follows in an analogous way.

For the last assertion, fixed $x\in Q$ take $y_n\in\partial^cv_n(x)$ such that $c(x,\cdot)$ reaches on $\partial^cv_n(x)$ its minimal value $c(x,y_n)$. By definition of $\partial^cv_n$ we have
$$v_n(x)+v_n^c(y_n)=c(x,y_n)$$
and we may assume that $y_n\to y$ in $Q$. By the first part of the lemma we may pass to the limit and deduce that
$$v(x)+v^c(y)=c(x,y)$$
which gives $y\in\partial^cv(x)$ and
$$\min_{\partial^cv(x)}c(x,\cdot)\le c(x,y)=\liminf_n c(x,y_n)=\liminf_n\Big(\min_{\partial^cv_n(x)}c(x,\cdot)\Big)$$
as required.
\end{proof}

We reformulate now problem \eqref{maxp} by considering the functional
\begin{equation}\label{functi}
I(v)=\int_Q\Big[v(x)-\min_{y\in\partial^cv(x)}c(x,y)\Big]\,df(x)
\end{equation}
on the admissible class
$$\B=\big\{v\hbox{ $c$-concave, }0\le v(x)\le v_0(x)\big\}$$
where
$$v_0(x)=\inf\big\{c(x,y)+p_0(y)\ :\ y\in Q\big\}.$$
By Lemma \ref{unifcont} the class $\B$ is compact for the uniform convergence, and by Lemma \ref{vc} the optimization problem
\begin{equation}\label{maxi}
\max\big\{I(v)\ :\ v\in\B\big\}
\end{equation}
admits a solution $v_{opt}$.

We can now come back to the initial problem \eqref{maxp} and deduce that it admits an optimal solution $p_{opt}$. Indeed, if $v \in \B$, then $p:=-v^c\in \A$ (since   for any function $u$ we have $(u^c)^c\ge u$) and $I(v)=F(p)$. Moreover, it is easy to check that $T_p(x)\subset \partial^c v_p(x)$ for every $p\in\A$. Thus,  $p_{opt}:=-v^c_{opt}$ actually solves \eqref{maxp} since for every $p\in\A$ one has
$$F(p)\leq \int_Q\Big[v_p(x)-\min_{y\in\partial^c v_p(x)}c(x,y)\Big]\,df(x)=I(v_p)\le I(v_{opt})=F(p_{opt}).$$

\section{Examples}\label{secex1}

So far, we have been in a rather abstract framework and it is time now to look at some special cases where the problem takes a more tractable form. 

\subsection{The case cost equal to distance}\label{distance}

We consider here the particular case when the cost function $c(x,y)$ is given by a distance $d(x,y)$ on $Q$; we shall see that in this situation the solution $p_{opt}$ above can be recovered in an explicit way. We denote by $\lip_{1,d}(Q)$ the class of all Lipschitz functions in $Q$ for the distance $d$ whose Lipschitz constant does not exceed 1.

\begin{theorem}\label{lipone}
In the case $c(x,y)=d(x,y)$ the optimal solution is given by
$$p_{opt}(x)=\max\big\{p(x)\ :\ p\in\lip_{1,d}(Q),\ p\le p_0\big\}.$$
\end{theorem}

\begin{proof}
We first notice that in this case the class of $c$-concave functions coincides with the class $\lip_{1,d}(Q)$. Moreover, as we have seen in the reduction from problem \eqref{maxp} to problem \eqref{maxi}, we may limit ourselves to consider only functions which are of the form $-v^c$ where $v$ is $c$-concave. In our case this allows us to limit the class of admissible $p$ to $\lip_{1,d}(Q)$.

Due to the tie-breaking rule it is easy to see that for $p\in\lip_{1,d}(Q)$ it is $T_p(x)=x$, which gives to the cost functional $F$ the simpler form
$$F(p)=\int_Q p(x)\,df(x).$$
Maximizing the previous expression in the class of functions in $\lip_{1,d}(Q)$ which are bounded by $p_0$ provides the solution
$$p_{opt}(x)=\max\big\{p(x)\ :\ p\in\lip_{1,d}(Q),\ p\le p_0\big\}$$
as required.
\end{proof}

\begin{remark}
We notice that in the case $c(x,y)=d(x,y)$ above the optimal pricing pattern $p_{opt}$ does not depend on the distribution $f$ of customers. Note also the explicit formula for  the optimal price:
$$p_{opt}(x)=\inf\{p_0(y)+d(x,y)\ :\ y\in Q\}\qquad\forall x\in Q.$$
\end{remark}

\begin{remark}
When $Q$ is a subset of  the Euclidean space $\R^N$, then Theorem \ref{lipone} in particular applies to the concave case where $d(x,y)=|x-y|^{\alpha}$ with $\alpha\in (0,1]$ since such costs are in fact metrics.
\end{remark}

\subsection{The case of a strictly convex cost}\label{strconvex}

We consider now the case  $Q:=\overline{\Omega}$ where $\Omega$ is some open bounded subset of the Euclidean space $\R^N$ and $c(x,y)=h(x-y)$ where $h$ is a nonnegative smooth and strictly convex function. In this framework, a $c$-concave function $v$ can be represented as:
\begin{equation}\label{envcconv}
v(x):=\min \{ h(x-y)-v^c(y),\; y\in Q\},\; \forall x\in Q.
\end{equation}
By the smoothness of $h$, the compactness of $Q$ and Lemma \ref{unifcont} ensure that $v$ is Lipschitz continuous on $Q$ hence Lebesgue a.e. differentiable on $\Omega$ by Rademacher's theorem. For every point $x\in\Omega$ of differentiability of $v$ and every $y\in \partial^c v(x)$, it is easy to check that from \pref{envcconv} one has:
\[\nabla v(x):=\nabla h(x-y)\] 
and since $h$ is strictly convex this can be rewritten as:
\begin{equation}\label{envelope}
y=x-\nabla h^*(\nabla v(x))
\end{equation}
where $h^*$ stands for the Legendre transform of $h$.  This proves that for every $c$-concave function $v$, $\partial^c v$ is in fact single-valued on a set of full Lebesgue measure. Now further assuming that $f$ is absolutely continuous with respect to the Lebesgue measure on $\Omega$, we can rewrite the profit functional in a more familiar form:
$$I(v)=\int_\Omega\left[v-h(\nabla h^*(\nabla v))\right]\,df
=\int_\Omega\left[v+h^*(\nabla v)-\nabla v\cdot\nabla h^*(\nabla v)\right]\,df.$$
If we further restrict our attention to the quadratic case, namely $c(x,y):=|x-y|^2/2$ and $\Omega$ is convex, it is easy to see that $v$ is $c$-concave on $\overline{\Omega}$ if and only if the function $w$ defined by
\[w(x):=\frac{1}{2} |x|^2-v(x),\qquad\forall x\in\Omega\]
is convex and satisfies 
\[ \nabla w(x)\in Q \mbox{ for a.e. } x\in \Omega.\]
Of course the constraint $v\leq v_0$ translates into $w\ge w_0$ with $w_0(x):= |x|^2/2-v_0(x)$. Putting everything together, we then see that $v$ solves \pref{maxi} if and only if $v(x)=|x|^2/2-w(x)$ and $w$ solves the following:
\begin{equation}\label{calcvarcvex}
\inf_{w\in\C} K(w)\mbox{ where }K(w):=\int_{\Omega}\left[\frac{1}{2}|\nabla w|^2-x\cdot\nabla w+w\right]\,df 
\end{equation}
and
$$\C:=\{w:\Omega\to\R,\ w\mbox{ convex,}\ w\ge w_0,\ \nabla w\in Q\mbox{ a.e.}\}.$$
Problems of the calculus of variations subject to a convexity constraint with a very similar structure as \pref{calcvarcvex} arise in the monopoly pricing model of Rochet and Chon\'e (\cite{rochetchone}).
Note also that by strict convexity, \pref{calcvarcvex} possesses a unique solution. 

\subsection{The quadratic case in dimension one}
We now consider problem \pref{calcvarcvex} in the special unidimensional case where  $\Omega=(0,1)$, $df=dx$ and $w_0\equiv 0$ (which corresponds to the price bound $p_0(x)=x-x^2/2$). The problem amounts to maximize $K(w)$ among convex, nondecreasing and $1$-Lipschitz functions $w$. It is obvious that one necessarily has $w(0)=0$ at the optimum, which setting $q:=w^{\prime}$ and integrating by parts enables us to write
\[K(w)=\int_0^1\left[\frac{1}{2} q(x)^2+(1-2x)q(x)\right]\,dx\] 
and the previous integral has to be minimized among nondecreasing functions $q(x)$ taking values in $[0,1]$. By a straightforward computation, the infimum is attained for $q_{opt}(x)=(2x-1)_+$, so that integrating we find $w_{opt}$ and then by $v_{opt}(x):=x^2/2-w_{opt}(x)$. Summarizing, we have obtained
$$v_{opt}(x)=\left\{\begin{array}{lll}
x^2/2  &\mbox{ if } &  x\in [0,1/2]\\
-x^2/2+x-1/4 & \mbox{ if } &  x\in [1/2,1].
\end{array}\right.$$
Finally, the optimal price is obtained by the formula $p_{opt}=-v_{opt}^c(x)$ which simply yields here $p_{opt}(x)=x/2-x^2/4=p_0(x)/2$.

\section{Problem formulation in the second case}\label{secondpb}

In what follows, $Q$ will be some compact metric space (the economic region), $Q_0$ is some open subset of $Q$ (the subregion where prices are fixed) and $p_0$ is a nonnegative l.s.c. function defined on $Q_0$ ($p_0$ is the fixed price system in $Q_0$). We are also given a transportation cost function $c$ assumed to be continuous and nonnegative on $Q\times Q$ and such that $c(x,x)=0$ for all $x\in Q$.  Finally, $f$ is a nonnegative Radon measure on $Q$ of positive mass that models the repartition of customers in $Q$. We set $Q_1:=Q\setminus Q_0$, this (compact) subregion being the one where prices have to be determined by the agent. 

\smallskip

The unknown of the problem is the pricing pattern $p$ in the following class:
$$\A:=\{p:Q\to\R,\ p=p_0\mbox{ on }Q_0,\ p\mbox{ l.s.c. on }Q\}.$$
Once a price $p$ has been fixed by the agent, consumers living at $x$ purchase the good so as to minimize their total expenditure i.e. price plus commuting cost, which leads to define, analogously to what done in Section \ref{firstpb},
\begin{equation}\label{defvptp2}
\left\{\hskip-0.4truecm\begin{array}{ll}
&v_p(x):=\min_{y\in Q}\big\{c(x,y)+p(y)\big\},\\
&T_p(x):=\big\{y\in Q\ :\ c(x,y)+p(y)=v_p(x)\big\}.
\end{array}\right.
\end{equation}
By our l.s.c. and compactness assumptions $T_p(x)$ is a nonempty compact subset of $Q$ and moreover the graph of $T_p$ is compact as the argmin of some l.s.c. function on $Q\times Q$. Note that $T_p(x)$ is the set of locations where consumers living at $x$ rationally choose to purchase the good. If $T_p(x)\subset Q_0$ then all the profit generated by the consumers of $x$ goes to the runner of region $Q_0$. We thus define:
\begin{equation}\label{defom}
\left\{\hskip-0.4truecm\begin{array}{ll}
&\Omega_0(p):=\big\{x\in Q\ :\ T_p(x)\subset Q_0\big\},\\
&\Omega_1(p):=\big\{x\in Q\ :\ T_p(x)\cap Q_1\ne\emptyset\big\}.
\end{array}\right.
\end{equation}
When $x\in\Omega_0(p)$, the agent makes no profit on consumers of $x$; when $x\in\Omega_1(p)$, we assume as \emph{tie-breaking rule} that the consumers living at $x$ go to a transportation-minimizing (or equivalently to a price-maximizing) location $y\in T_p(x)$:
$$y\in\argmin_{T_p(x)\cap Q_1} c(x,\cdot)=\argmax_{T_p(x)\cap Q_1} p(\cdot)$$
(notice that every $y\in T_p(x)$ yields the same minimal expenditure to $x$). 

With the previous notations, we see that the optimal pricing problem amounts to the maximization problem
\begin{equation}\label{pbmeenp}
\max\big\{\Pi(p)\ :\ p\in\A\big\}
\end{equation}
where
$$\Pi(p):=\int_{\Omega_1(p)}\left(\max_{y\in T_p(x)\cap Q_1}p(y)\right)\,df(x).$$
By the definition of $T_p$ and $v_p$, one has $v_p(x)=c(x,y)+p(y)$ for all $y\in T_p(x)$, hence the profit functional can be rewritten as: 
\[\Pi(p)=\int_{\Omega_1(p)} \left( v_p(x)-\min_{y\in T_p(x)\cap Q_1} c(x,y) \right) df(x).\]
Defining for all $x\in \Omega_1(p)$
\begin{equation*}
\left\{\hskip-0.4truecm\begin{array}{ll}
&H_p(x):=\max_{y\in T_p(x)\cap Q_1} p(y),\\
&G_p(x):= \min_{y\in T_p(x)\cap Q_1} c(x,y)
\end{array}\right.
\end{equation*}
we may then rewrite in a more synthetical way the profit as
$$\Pi(p)=\int_{\Omega_1(p)}H_p(x)\,df(x)
=\int_{\Omega_1(p)}(v_p(x)-G_p(x))\,df(x).$$

\begin{remark}\label{rk1}
So far, we have not assumed that $p$ has to be nonnegative, in fact this constraint is unnecessary since it will directly follow from the maximization problem \pref{pbmeenp}. Indeed if $p\in\A$ then $p_+:=\max(p,0)$ is also in $\A$ and $\Pi(p_+)\ge\Pi(p)$. If $H_p\le0$ on $\Omega_1(p)$, this claim is obvious. We may then assume that $\{H_p\ge0\}\cap\Omega_1(p)\ne\emptyset$. Let $x\in\Omega_1(p)$ be such that $H_p(x)\ge0$ and let $y\in T_p(x)\cap Q_1$ be such that $H_p(x)=p(y)=p_+(y)$; we have $v_p(x)=c(x,y)+p(y)=c(x,y)+p_+(y)\ge v_{p_+}(x)$ and since $v_p\le v_{p_+}$ this yields $v_{p_+}(x)=c(x,y)+p_+(y)$ which implies $x\in \Omega_1(p_+)$, $y\in T_{p_+}(x)$ and $H_{p_+}(x)\ge H_p(x)$. We then have
$$\Pi(p)\le\int_{\Omega_1(p)\cap\{H_p\ge0\}} H_p(x)\,df(x)
\le\int_{\Omega_1(p_+)} H_{p_+}(x)\,df(x)=\Pi(p_+).$$
\end{remark}

\section{The existence result}\label{sec2}

\subsection{Generalized concavity}\label{gc}
To prove the existence of a maximizer in \pref{pbmeenp}, we reformulate the problem in terms of $v_p$ rather than on the price $p$ which a priori does not have special properties. To do that, it is convenient to use some notions of generalized concavity that are natural in our context (as well as in the Monge-Kantorovich theory). Before introducing formal definitions let us remark that for $p\in \A$, one can rewrite $v_p$ (defined by \pref{defvptp2}) as:
\[v_p(x)= v_0(x)\wedge w_p(x)\]
(where $a\wedge b$ denotes the minimum of the two real numbers $a$ and $b$) with
\begin{equation}\label{defwp}
\left\{\hskip-0.4truecm\begin{array}{ll}
&v_0(x):=\inf_{y\in Q_0}\big\{c(x,y)+p_0(y)\big\},\\
&w_p(x):=\inf_{y\in Q_1}\big\{c(x,y)+p(y)\big\}.
\end{array}\right.
\end{equation}
From the previous formula we see that $w_p$ can be represented as the pointwise infimum of a family of functions $x\mapsto c(x,y)+p(y)$ where the parameter $y$ takes its values in $Q_1$. This suggests the following definition.

\begin{definition}
A function $w:Q\to\R$ is called $(Q_1,c)$-concave if there exists a function $u:Q_1\to\R$ bounded from above such that 
\begin{equation}\label{defccon}
w(x)=\inf_{y\in Q_1}\{c(x,y)-u(y)\},\; \forall x\in Q.
\end{equation}
\end{definition}

If $w$ is $(Q_1,c)$-concave there exists a kind of minimal representation (as for the usual Legendre-Fenchel transform) of $w$ in the form \pref{defccon}. Indeed, using the $c$-transform (see Definition \ref{cconc})
$$w^c(y):=\inf_{x\in Q}\{c(x,y)-w(x)\}\qquad\forall y\in Q_1,$$
one has
\begin{equation}\label{ctransf}
w(x)=\inf_{y\in Q_1}\{c(x,y)-w^c(y)\}\qquad\forall x\in Q.
\end{equation}
Indeed, on the one hand, the definition of $w^c$ yields $w(x)+w^c(y)\le c(x,y)$ for every $(x,y)\in Q\times Q_1$, hence:
$$w(x)\le\inf_{y\in Q_1}\{c(x,y)-w^c(y)\}.$$
On the other hand, using the representation \pref{defccon} yields $u\le w^c$ on $Q_1$ hence
$$w(x)=\inf_{y\in Q_1}\{c(x,y)-u(y)\}\ge\inf_{y\in Q_1}\{c(x,y)-w^c(y)\}.$$
Analogously to what was done in Section \ref{firstpb}, for every $(Q_1,c)$-concave function $w$, the $(Q_1,c)$-superdifferential of $w$ at $x\in Q$ (denoted $\partial^{1,c} w(x)$) is defined by
\[\partial^{1,c} w(x):=\{y\in Q_1 \; : \; w(x)+w^c(y)=c(x,y)\}.\]
Since $Q$ and $Q_1$ are compact and $(Q_1,c)$-concave functions and their $c$-transforms are continuous, it is easy to see that for every $(Q_1,c)$-concave function $w$ and every $x\in Q$, $\partial^{1,c} w(x)$ is a nonempty compact subset of $Q_1$ and that $\{(x,y)\in Q\times Q_1 \; : \; y\in \partial^{1,c} w(x)\}$ is compact.

\subsection{Reformulation}

The aim of this subsection is to reformulate the maximization problem \pref{pbmeenp} in terms of $w=w_p$ only. Let $p\in \A$ be nonnegative (which is not restrictive in view of Remark \ref{rk1}) and write 
\[v_p:=v_0\wedge w\]
with $w=w_p$ and $v_0$ defined by \pref{defwp}. Then, let us define
\begin{equation}\label{deftilu}
\tilu(y):=\inf_{x\in Q}\{c(x,y)-w(x)\}\qquad\forall y\in Q_1;
\end{equation}
as already noticed, since $w=w_p$ is $(Q_1,c)$-concave we have
$$w(x):= \inf_{y\in Q_1} \{c(x,y)-\tilu(y)\}\qquad\forall x\in Q.$$
Now let us define
\[\tilp(y):=\left\{\begin{array}{lll}
p_0(y) &\mbox{ if } &  y\in Q_0\\
-\tilu(y) & \mbox{ if } &  y\in Q_1.
\end{array}
\right.\]
By construction $w_{\tilp}=w_p$ hence $v_{\tilp} =v_p$. The next proposition expresses that the profit is improved when one replaces $p$ by $\tilp$. This allows us to restrict the analysis to prices that are (up to a minus sign) a $c$-transform on the free region $Q_1$ and will enable us to express the problem in terms of $w$ only. More precisely, we have the following.

\begin{proposition}\label{reformul}
Let $p\in\A$, $p\ge0$ and let $\tilp$ and $w$ be defined  as above. Then one has
\begin{eqnarray}
&v_{\tilp} =v_p,\quad \tilp\le p\mbox{ on }Q_1,\quad \tilp\ge0\mbox{ on }Q,\label{prop1re}\\
&T_p(x)\cap Q_1\subset T_{\tilp}(x) \cap Q_1,\quad\forall x\in \Omega_1(p),\label{prop2re}\\
&\Omega_1(p) \subset \Omega_1(\tilp)=\{w\leq v_0\}, \label{prop3re}\\
&T_{\tilp}(x)\cap Q_1=\partial^{1,c} w(x),\quad\forall x\in\Omega_1(\tilp), \label{prop4re}
\end{eqnarray}
which imply 
\begin{equation}\label{pitamel}
\Pi(\tilp)\ge\Pi(p)
\end{equation}
and 
\begin{equation}\label{paw}
\Pi(\tilp)=\int_{\{w\leq v_0\}}\left(w(x)-\min_{y\in\partial^{1,c} w(x)} c(x,y)\right)\,df(x).
\end{equation}
\end{proposition}

\begin{proof}
We already know that $v_{\tilp} =v_p$. Using Subsection \ref{gc} we know that $u:=-p\leq \tilu=-\tilp$ on $Q_1$. Since $p\ge0$ and $w\ge0$ we have 
\[\tilu(y)= \inf_{x\in Q} \{c(x,y)-w(x)\}\leq \inf_{x\in Q} \{c(x,y)\}=0\]
which proves $\tilp\geq 0$ and \pref{prop1re}.

Now, fix $x\in \Omega_1(p)$ and $y\in T_p(x)\cap Q_1$; then $v_p(x)=v_{\tilp}(x)=c(x,y)+p(y)\ge w(x)$ which, since $p(y)\geq \tilp(y)$, implies that
$$w(x)=v_p(x)\ge c(x,y)+\tilp(y)\ge v_{\tilp}(x)$$
so that $y\in T_{\tilp}(x)$ which proves \pref{prop2re}. The previous argument also proves that $\Omega_1(p)\subset\Omega_1(\tilp)$. The fact that $\Omega_1(\tilp)=\{w\le v_0\}$ and \pref{prop4re} are obvious. Thanks to the fact that $\tilp\ge0$, the integrand in $\Pi(\tilp)$ is nonnegative and thanks to \pref{prop2re} we obtain $G_{\tilp}\leq G_p$ on $\Omega_1(p)$. By \pref{prop3re} we then have:
\[\Pi(\tilp)\geq \int_{\Omega_1(p)} (v_{\tilp}-G_{\tilp})df\geq \int_{\Omega_1(p) }(v_p-G_p)df=\Pi(p).\]
Finally, thanks to \pref{prop3re} and \pref{prop4re}, $\Pi(\tilp)$ can be rewritten as a function of $w$ only as in \pref{paw}.
\end{proof}

Proposition \ref{reformul} thus enables us to reformulate the initial problem \pref{pbmeenp} as:
\begin{equation}\label{pbmeenw}
\sup_{w\in \W} J(w)=\int_{\{w\leq v_0\}} \left(w(x)-\min_{y\in \partial^{1,c} w(x)} c(x,y)\right)\,df(x)
\end{equation}
where $\W$ is the set of all $(Q_1,c)$-concave functions. More precisely, if $w$ solves \pref{pbmeenw} then $\tilp$ defined by
\[\tilp(y):=\left\{\begin{array}{lll}
p_0(y) &\mbox{ if } &  y\in Q_0\\
-\tilu(y) & \mbox{ if } &  y\in Q_1.
\end{array}
\right.\]
solves \pref{pbmeenp}.

Notice that if $w\in\W$ then one has for every $(x_1,x_2)\in Q\times Q$
\begin{equation}\label{equic}
|w(x_1)-w(x_2)|\le\max_{y\in Q_1}|c(x_1,y)-c(x_2,y)|,
\end{equation}
which proves that $\W$ is an equicontinuous family. 

\subsection{Existence}

With the reformulation \pref{pbmeenw} (and the equicontinuity estimate \pref{equic}) at hand, we easily deduce the following existence result.

\begin{theorem}\label{exist}
Problem \pref{pbmeenw} admits at least one solution (hence so does problem \pref{pbmeenp}).
\end{theorem}

\begin{proof}
Let $(w_n)_n$ be some maximizing sequence of \pref{pbmeenw}; without loss of generality we may assume that the integrand in the definition of $J(w_n)$ is nonnegative (see Remark \ref{rk1}) and that $\min_Q (w_n- v_0)\le0$. By \pref{equic} we deduce that $(w_n)$ is uniformly bounded and equicontinuous. Thanks to Ascoli-Arzela's theorem, passing to a subsequence if necessary we may assume that $w_n$ converges uniformly to some $w$ which is easily seen to be $(Q_1,c)$-concave too. To prove that $w$ solves \pref{pbmeenw}, we first use Fatou's lemma:
\[\limsup J(w_n)\leq \int_Q \limsup \chi_{\{w_n\leq v_0\}} \left(w_n(x)-\min_{y\in \partial^{1,c} w_n(x)} c(x,y) \right)\,df(x).\]
It is therefore enough to prove that for every $x\in Q$

\begin{equation}\label{scsx}
\begin{array}{lll}
&\limsup\chi_{\{w_n\leq v_0\}}(w_n(x)-\min_{y\in\partial^{1,c} w_n(x)} c(x,y))\\
&\qquad\qquad\le\chi_{\{w\leq v_0\}}(w(x)- \min_{y\in \partial^{1,c} w(x)} c(x,y)).
\end{array}
\end{equation}
If $w(x)>v_0(x)$ the right-hand side vanishes and, since $w_n\to w$ uniformly, we have $w_n(x)>v_0(x)$ for $n$ large enough, so that the left-hand side vanishes too. Assume now that $w(x)\leq v_0(x)$, and let $y_n\in \partial^{1,c} w_n(x)$ be such that
\[c(x,y_n)=\min_{y\in \partial^{1,c} w_n(x)} c(x,y);\]
passing to a subsequence if necessary we may assume that $y_n$ converges to some $y\in \partial^{1,c} w(x)$, hence \pref{scsx} holds.
\end{proof}

\section{Examples}

\subsection{The eikonal case}\label{concave}

In this subsection, we investigate the particular case where $Q=\Omb$, the closure of a bounded open convex subset of $\R^d$ and the cost is the euclidean distance $c(x,y)=|x-y|$. As before we assume that $Q_0$ is an open subset of $Q$ and $Q_1=Q\setminus Q_0$.  As already noticed, in this case, the $c$-concave functions are simply the $1$-Lipschitz ones. As for the $(Q_1,c)$-concave ones,  it is easy to see that  $w$ is $(Q_1,c)$-concave if and only if it is $1$-Lipschitz on $Q$ and
\begin{equation}\label{repwp}
w(x)= \min_{y\in Q_1}\{|x-y|+w(y)\},\qquad\forall x\in Q.
\end{equation}
Now, let  $x\in Q_0$ be a point of differentiability of $w$ and let $y\in Q_1$ (so that $x\neq y$) be such that $w(x)=|x-y|+w(y)$ (i.e. $y\in \partial^{1,c} w(x)$), then one has
\begin{equation}
\nabla w(x)=\frac{x-y}{|x-y|}\mbox{ and  there exists }\lambda>0\mbox{ such that }x-\lambda\nabla w(x)\in Q_1
\end{equation}
so that 
\begin{equation}\label{condsurw}
|\nabla w(x)|=1\mbox{ and  }\nabla w(x)\in\R_+(x- Q_1).
\end{equation}
By Rademacher's Theorem, \pref{condsurw} holds a.e. on $Q_0$. In particular $w$ is an a.e. solution of the eikonal equation $|\nabla w|=1$ on $Q_0$. Let $x\in Q_0$ be a point of differentiability of $w$, $y\in \partial^{1,c} w(x)$ and $\lambda=|x-y|$, then with the fact that $w$ is $1$-Lipschitz, it is easy to check that $w(x)-w(x-t\nabla w(x))=t$, for every $t\in [0,\lambda]$ (i.e. $w$ grows at the maximal rate $1$ on the segment $[x-\lambda \nabla w(x),x]$). In particular,  choosing $t\in[0,\lambda]$ such that $x-t\nabla w(x)\in \partial Q_0$ yields:
$$w(x)\ge\min_{y\in \partial Q_0}\{|x-y|+w(y)\}.$$
By density, this inequality actually holds for all $x\in Q_0$, and the converse inequality follows immediately from \pref{repwp}.We thus have proved that if   $w$ is $(Q_1,c)$-concave then
\begin{equation}\label{repwpb}
w(x)= \min_{y\in \partial Q_0}\{|x-y|+w(y)\},\qquad\forall x\in Q_0.
\end{equation}
It is well-known (see \cite{barles}) that \pref{repwpb} implies that $w$ is  a \emph{viscosity} solution of the eikonal equation on $Q_0$.  Now, conversely, assume that $w$ is $1$-Lipschitz on $Q$ and a viscosity solution of the eikonal equation on $Q_0$ and define
\begin{equation}\label{defu}
u(x)= \min_{y\in Q_1}\{|x-y|+w(y)\},\qquad\forall x\in Q.
\end{equation}
then $u=w$ on $Q_1$ (in particular on $\partial Q_0$) and by the same argument as above $u$ is a viscosity solution of the eikonal equation on $Q_0$. A standard comparison argument (e.g. Theorem 2.7 in \cite{barles}) yields $u=w$ on $Q_0$ so that $w$ is $(Q_1,c)$-concave. This proves that the set of $(Q_1,c)$ concave functions is:
\begin{equation}\label{formW}
\W=\{w:Q\to\R,\ \mbox{$w$ $1$-Lipschitz on $Q$ and $|\nabla w|=1$ on $Q_0$}\}
\end{equation}
where the eikonal equation has to be understood in the viscosity sense. Let us also remark that the condition $\nabla w(x)\in\R_+(x-Q_1)$ a.e. in $Q_0$ is in fact hidden in the definition of a viscosity solution (equivalently in formula \pref{repwpb}).

Getting back to our optimization problem \pref{pbmeenw}, it is natural to introduce for every $x\in Q$ and $\nu\in S^{d-1}$ (the unit sphere of $\R^d$) the quantity:
\begin{equation}
\lambda(x,\nu):=\inf\{\lambda\ge0\ :\ x-\lambda\nu\in Q_1\}.
\end{equation}
For $w\in\W$, we then have for a.e. $x\in Q$
$$\min_{y\in\partial^{1,c}w(x)}|x-y|=\lambda(x,\nabla w(x))$$
Assuming that $f$ is absolutely continuous with respect to the Lebesgue measure on $\Omega$ and defining $v_0$ by \pref{defwp}, for $w\in \W$, the profit functional $J$ is then given by:
$$J(w):=\int_{\{w\le v_0\}}\Big(w(x)-\lambda(x,\nabla w(x))\Big)\,df(x)$$
which has to be maximized over $\W$ defined by \pref{formW}. Now, our aim is to transform the previous problem in terms of the values of $w$ on $\partial Q_0$ only. Of course, because of \pref{repwpb}, the behavior of $w$ on $Q_0$ is fully determined by its trace on $\partial Q_0$. In order to treat the behavior on $Q_1$, we need the following result.

\begin{lemma}
Let $w\in\W$ and define
$$u(x):=\min_{y\in Q_1}\{|x-y|+u(y)\wedge v_0(y)\},\qquad\forall x\in Q,$$
then $u\in \W$ and $J(u)\ge J(w)$. 
\end{lemma}

\begin{proof}
Obviously $u\in \W$ and $u=u
\wedge v_0$ on $Q_1$ hence the integrand in the definition of $J$ is larger on $Q_1$ for $u$ than for $w$ (recall that $v_0\ge0$). If $x\in Q_0$ is such that $w(x)>v_0(x)$, then the same conclusion holds. Now, if $x\in Q_0$  is such that $w(x)\le v_0(x)$, then we write $u(x)=u(y)+|x-y|$ with $y\in\partial^{1,c} u(x)$, if $w(y)\le v_0(y)$ then  $u(x)=w(y)+|x-y|\ge w(x)$ and if $w(y)\ge v_0(y)$ then $u(x)=v_0(y)+|x-y|\ge v_0(x)\ge w(x)$. Since $u\le w$, in both cases we then have $u(x)=w(x)$ which proves that $Q_0\cap \{w\leq v_0\}\subset Q_0\cap \{u=w\}$. In particular, $u-\lambda(x,\nabla u)=w-\lambda(x,\nabla w)$ a.e. on $Q_0\cap \{w\leq v_0\}$ which proves the desired result. 
\end{proof}

Let $w\in \W$ and let $\phi$ be the trace of $w$ on $\partial Q_0$, thanks to the previous Lemma we may assume that $w\leq v_0$ on $Q_1$ so that:
$$J(w)=\int_{Q_1} w\,df +\int_{Q_0\cap\{w\le v_0\}}\Big(w(x)-\lambda(x,\nabla w(x))\Big)\,df(x).$$
Because of \pref{repwpb}, the second term only depends on $\phi$, and the first one is monotone in $w$ hence for a given $\phi$ ($1$-Lipschitz and smaller than $v_0$) it is maximized by the largest $1$-Lipschitz function on $Q_1$ which has $\phi$ as trace on $Q_0$ and is below $v_0$ i.e. simply
$$w(x)=\min_{y\in\partial Q_0} \{|x-y|+\phi(y)\},\qquad\forall x\in Q_1.$$
Since the previous formula also holds for $x\in Q_0$ by  \pref{repwpb}, we define for every $1$-Lipschitz function $\phi$ on $\partial Q_0$ such that $\phi\le v_0$ the state equation
\begin{equation}\label{state}
w_\phi(x):=\min_{y\in\partial Q_0} \{|x-y|+\phi(y)\},\qquad\forall x\in Q.
\end{equation}
The profit maximization \pref{pbmeenw} can thus be reformulated as the following nonstandard optimal control problem where the control is the price $\phi$ on the interface $\partial Q_0$:
\begin{equation}\label{controlform}
\sup_{\phi\in\Phi} J(w_\phi)=\int_{Q_1} w_\phi\,df+\int_{Q_0\cap\{w_\phi\le v_0\}}\Big(w_\phi(x)-\lambda(x,\nabla w_\phi(x))\Big)\,df(x)
\end{equation}
where the class of admissible boundary controls $\Phi$ consists of all $1$-Lipschitz functions $\phi$ on $\partial Q_0$ such that $\phi\leq v_0$ and the state equation is \pref{state}.

For example if $Q$ is the unit ball of $\R^d$ and $Q_1$ its boundary, then the maximization problem \pref{pbmeenw} becomes maximizing:
$$J(w):= \int_{\{w\le v_0\}}\Big(w(x)-x\cdot\nabla w(x)-\sqrt{(x\cdot\nabla w(x))^2+|x|^2-1}\Big)\,df(x)$$
in the set of viscosity solutions of the eikonal equation $|\nabla w|=1$ on the unit ball.  Note that this is a highly nonconvex variational problem, which as previously may be reformulated as maximizing $J(w_{\phi})$ among $1$-Lipschitz functions $\phi$ on $\partial Q_0$ such that $\phi\le v_0$.

\subsection{The one-dimensional case}\label{onedim} 

In the one dimensional case, the eikonal equation has a very simple structure which makes problem \pref{controlform} much simpler. In particular, if $\partial Q_0$ is finite then the maximization of \pref{controlform} reduces to a finite dimensional problem, since the control in this case is simply given by the values of $w$ on the finite set $\partial Q_0$. For instance let us take $Q=[0,1]$, $Q_0=(\alpha,\beta)$ with $0\le\alpha<\beta\le1$. For simplicity let us also assume that $p_0$ is constant on $Q_0$ and that $f$ is a probability that does not charge points. Then the solutions of \pref{controlform} only depend on the two scalars $p_1:=w(\alpha)$ and $p_2:=w(\beta)$ subject to the constraints:
\begin{equation}\label{contdf} 
p_1\le p_0,\quad p_2\le p_0,\quad |p_2-p_1|\le\beta-\alpha.
\end{equation}
For such a control $(p_1,p_2)$ the function $w_{(p_1,p_2)}$ has the following W-like shape:

\centerline{\includegraphics[width=2.5in, angle=270]{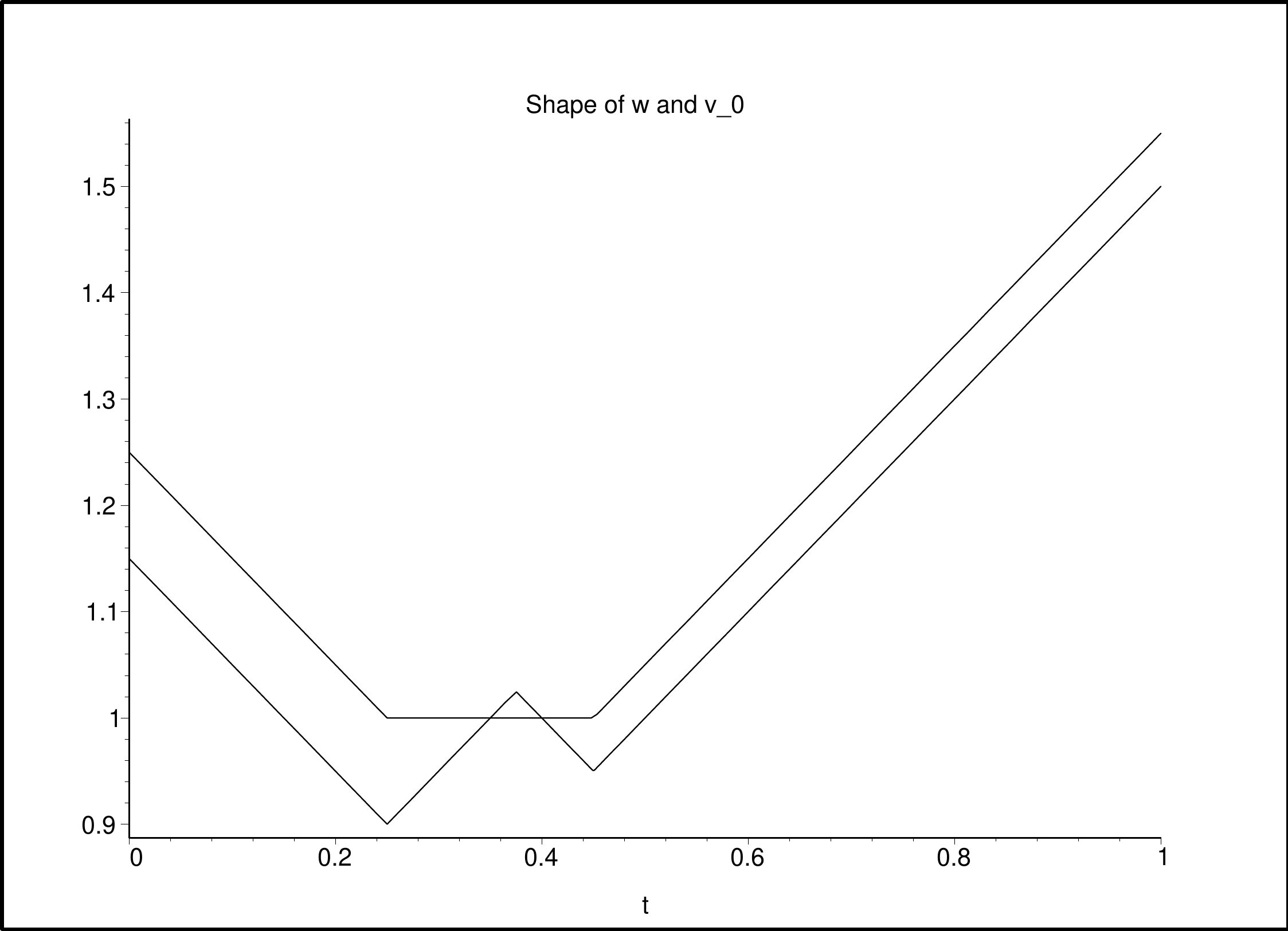}}

\bigskip

The function $\lambda(x,\nu)$ is in this case
$$\lambda(x,\nu)=\left\{
\begin{array}{ll}
x-\alpha&\mbox{if $x\in]\alpha,\beta[$ and $\nu=1$,}\\
\beta-x&\mbox{if $x\in]\alpha,\beta[$ and $\nu=-1$,}\\
0&\mbox{otherwise,}
\end{array}\right.$$
and the corresponding profit can be explicitly computed as a function of $(p_1,p_2)$:
\[\begin{split}
\int_0^\alpha (p_1+\alpha-s)\,df(s)+\int_{\alpha}^{(s_0\wedge s_1)(p_1,p_2)} p_1 \,df(s)\\
+\int_{(s_0\vee s_2)(p_1,p_2)}^\beta p_2\,df(s)+\int_\beta^1 (p_2+s-\beta)\,df(s)
\end{split}\]
where
$$\left\{
\begin{array}{l}
s_1(p_1,p_2)=p_0-p_1+\alpha,\\
s_2(p_1,p_2)=p_2-p_0+\beta,\\
s_0(p_1,p_2)=\frac{1}{2}(p_2-p_1+\beta+\alpha).
\end{array}
\right.$$
Defining $F$ the cumulative function of $f$ (i.e. $F(t)=f([0,t])$), solving \pref{controlform} then amounts to maximize:
$$p_1F((s_0\wedge s_1)(p_1,p_2))+p_2\Big(1-F((s_0\vee s_2)(p_1,p_2))\Big)$$
subject to the constraints \pref{contdf}. For example, if $\alpha=0$, $\beta=1$ (i.e. the price $p(x)$ has to be chosen only at the boundary of $Q$) and $f$ is uniform, then there is a unique optimal strategy given by $p_1=p_2=\frac{p_0}{2}\vee (p_0-\frac{1}{2})$.

\section{Concluding remarks and related problems}\label{lastsec}

In this section we propose some further developments of the optimization problems above that could be investigated. It is not our goal to enter into the details, which could possibly be treated in a future paper.

The model problems considered in the previous sections could also be used to describe a two (or more) players game, where each player operates only on its own region and considers the prices on the other regions as fixed. More precisely, assume that $Q=A\cup B$ where $A$ and $B$ are two compact sets with no interior point in common (although this is not essential for what follows). On $A$ and $B$ two agents (for instance the central governments of two different countries) operate and initially two price functions $p_0(x)$ and $q_0(x)$ are present on $A$ and $B$ respectively.

At a first step the agent that operates on $A$ modifies its price on $A$ considering $q_0$ fixed on $B$ and maximizes its income choosing an optimal price function $p_1$; then the agent that operates on $B$ plays its move considering $p_1$ fixed on $A$ and maximizing its income through the choice of an optimal price function $q_1$. The game continues in this way then providing price functions $p_n$ and $q_n$ defined on $A$ and $B$ respectively.

An interesting issue would be the study of the convergence of the sequences $(p_n)$ and $(q_n)$ to price strategies $p$ and $q$ that the two agents do not have the interest to modify any more.

A related alternative is to consider the competitive problem between the agents operating on $A$ and $B$ as a two-persons game (see for instance \cite{ft} or \cite{osborne}), which is not zero-sum, where the strategy of each player is the pricing function on the region he controls.  One has to be cautious in precisely defining the payoff functions when some customers are indifferent between being the good in $A$ or in $B$. In such a case, one can for instance impose, as tie-breaking rule, that each customer shops in his own region, and for simplicity we assume $f(A\cap B)=0$.

For respective price strategies $p$ (prices on $A$) and $q$ (prices on $B$), define for all $x\in Q$
\[\begin{split}
&v_p(x):=\inf_{y\in A} \{c(x,y)+p(y)\},\\
&w_q(x):=\inf_{z\in B} \{c(x,z)+q(z)\},\\
&T_p(x):=\{y\in A\;  : \; v_p(x)=c(x,y)+p(y)\},\\  
& S_q(x):=\{z\in B\;  : \; w_q(x)=c(x,z)+q(z)\}.
\end{split}\]
Under our tie-breaking rule, the payoff functions for the two players are then given by
\[\begin{split}
&\Pi_A(p,q):=\int_{\{v_p<w_q\}}\Big(\max_{y\in T_p(x)} p(y)\Big)\,df(x)
+\int_{\{v_p=w_q\}\cap A}\Big(\max_{y\in T_p(x)} p(y)\Big)\,df(x)\\
&\Pi_B(p,q):=\int_{\{w_q<v_p\}}\Big(\max_{z\in S_q(x)} q(z)\Big)\,df(x)
+\int_{\{v_p=w_q\}\cap B}\Big(\max_{z\in S_q(x)} q(z)\Big)\,df(x).  
\end{split}
\]
Defining admissible strategies as pairs of nonnegative and l.s.c. functions on $A$ and $B$ (possibly also satisfying additional constraints), an interesting issue is then to find Nash equilibria (see for instance \cite{ft} or \cite{osborne}) for the payoffs $(\Pi_A,\Pi_B)$, that is a pair of admissible strategies $p^*$ and $q^*$ such that
\[\Pi_A(p^*,q^*)\geq \Pi_A(p,q^*), \;  \Pi_B(p^*,q^*)\geq \Pi_B(p^*,q), \; \forall \mbox{  admissible strategies $p$ and $q$} . \]  
This is a priori a complicated problem because Kakutani's fixed-point Theorem does not apply here because of the tie-breaking rule which induces discontinuities. Of course, one can extend the framework to more than two-players, introduce mixed-strategies... The analysis of spatial competition is an important issue in economics since Hotelling's celebrated model  \cite{hotelling}  and one may relate the equilibrium problem described above to this line of research. The study of the general Nash problem is left for future research, but we give an elementary example where the solution is very simple and intuitive.

\begin{example}
Let $Q=[0,1]$, $A=[0,1/2]$, $B=[1/2,1]$, $\alpha\in (0,1)$ and $c(x,y)=|x-y|$. As explained in Section \ref{concave}, given the strategy of the second (respectively first) player, the first  (resp. second) one  optimally choses a pricing function of the form $a+1/2-x$ (resp. $b+x-1/2$) for $x\in A$ (resp. for $x\in B$). At a Nash equilibrium one must have $a=b$ (if $a>b$ then $A$ makes zero profit as well as $B$ makes zero profit if $a<b$). Finally, the common value $a=b$ has to be $0$, since if for instance $a>0$ then $B$ can charge a slightly lower price $a-\eps$ at the border point $1/2$ then getting the whole demand and increasing his profit for $\eps$ small enough. In this simple case there is then a unique Nash equilibrium $p(x)=1/2-x$ and $q(x)=x-1/2$, no matter what the population distribution is. The equilibrium price is plotted in the next figure.

\centerline{\includegraphics[width=2.5in, angle=270]{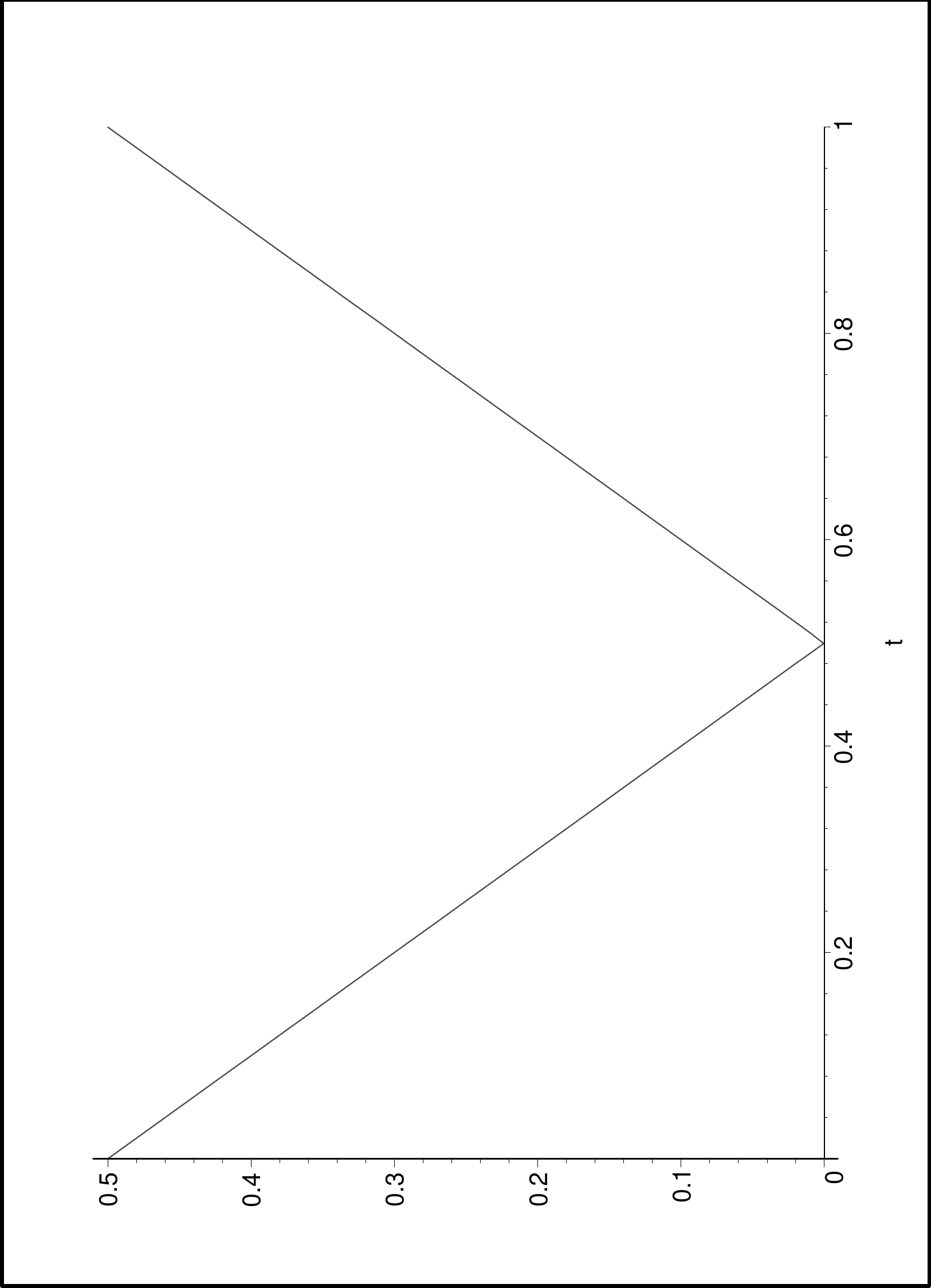}}

\end{example}

%

\bigskip
{\small
\begin{minipage}[t]{6.3cm}
Giuseppe Buttazzo\\
Dipartimento di Matematica\\
Universit\`a di Pisa\\
Largo B. Pontecorvo, 5\\
56127 Pisa - ITALY\\
{\tt buttazzo@dm.unipi.it}
\end{minipage}
\begin{minipage}[t]{6.3cm}
Guillaume Carlier\\
CEREMADE\\
Universit\'e de Paris-Dauphine\\
Place du Mar\'echal De Lattre De Tassigny\\
75775 Paris Cedex 16 - FRANCE\\
{\tt carlier@ceremade.dauphine.fr}
\end{minipage}
}

\end{document}